\documentclass[11pt,a4paper]{article}
\usepackage[utf8]{inputenc}
\usepackage[T1]{fontenc}
\usepackage[english]{babel}
\usepackage{amsmath,amssymb,amsthm,mathtools}
\usepackage[a4paper,margin=1in]{geometry}
\usepackage{parskip}
\usepackage{tikz}
\usetikzlibrary{positioning,arrows.meta,calc,decorations.pathmorphing}
\usepackage{microtype}
\usepackage[colorlinks=true,citecolor=blue,linkcolor=blue,urlcolor=blue]{hyperref}

\theoremstyle{plain}
\newtheorem{theorem}{Theorem}[section]
\newtheorem{proposition}[theorem]{Proposition}
\newtheorem{lemma}[theorem]{Lemma}

\theoremstyle{definition}
\newtheorem{definition}[theorem]{Definition}
\newtheorem{example}[theorem]{Example}
\theoremstyle{remark}
\newtheorem{remark}[theorem]{Remark}

\newcommand{\RR}{\mathbb{R}}
\newcommand{\NN}{\mathbb{N}}
\newcommand{\drho}{d_{\rho}}
\newcommand{\lerho}{\leq_{\rho}}

\title{A Kleene-type PA fixed point on partial metrics}

\author{%
  D.J.~Leko\thanks{AI Research and Innovation Nexus for Africa (AIRINA Labs), B\'enin; African Center for Advanced Studies (ACAS), P.O. Box 4477, Yaound\'e, Cameroon.}%
  \and Y.U.~Gaba\thanks{AI Research and Innovation Nexus for Africa (AIRINA Labs), B\'enin; African Center for Advanced Studies (ACAS), P.O. Box 4477, Yaound\'e, Cameroon.}%
}

\date{November 2026}

\begin{document}
\maketitle

\begin{abstract}
\noindent In a Matthews partial metric space, the self-distance
$\rho(x, x)$ need not be zero; the induced metric $\drho$ obtained by
zeroing the self-distances of $\rho$ turns $(X, \rho)$ into a standard
metric space on which most generalized contraction fixed-point theorems
reduce to their classical counterparts. The order-theoretic
Matthews--Kleene theorem for operators monotone with respect to the
canonical order $\lerho$ attached to $\rho$ sits outside the perimeter
of this reduction: its conclusion names a $\lerho$-supremum, and
$\lerho$ is erased by $\drho$. We formulate a path-averaged (PA)
condition in the sense of Fabiano on self-distances along a single
$\lerho$-monotone orbit, strictly weaker than the Matthews contraction,
and prove that under PA the supremum $x^{\ast} =
\sup_{\lerho}\{f^{n}(x_{0})\}$ satisfies $\rho(x^{\ast}, x^{\ast}) = 0$.
The role of PA in the theorem is narrow: existence, the supremum
characterization, and minimality above $x_{0}$ follow from Kleene's
theorem on a chain-complete poset without PA; PA is used only to force
the self-distance sequence to zero and so to secure the
strong-fixed-point property. A worked example on the lattice of bounded
intervals and an off-chain counterexample on which the Matthews
contraction fails but PA succeeds illustrate the scope.

\vspace{0.3em}
\noindent \textit{Keywords:} partial metric, PA contraction, Kleene
fixed point, canonical order, reduction-resistance.

\vspace{0.3em}
\noindent \textit{MSC 2020:} 47H10, 54E99, 06B23.
\end{abstract}

\section{Introduction}
\label{sec:intro}

Matthews partial metrics~\cite{matthews1994} were introduced for the
denotational semantics of programming languages. The self-distance
$\rho(x, x)$ measures how far the information attached to $x$ is from
being totally defined, is not necessarily zero, and vanishes only at
points that carry complete information. The framework has produced in
thirty years a wide literature of fixed-point results~\cite{romaguera2010,
oltra2004, kirkShahzad2014}.

Kirk and Shahzad~\cite[\S14.6]{kirkShahzad2014} made the following
structural observation. Every partial metric $\rho$ on $X$ induces
a metric
\begin{equation}\label{eq:drho}
\drho(x, y) =
\begin{cases}
0 & \text{if } x = y, \\
\rho(x, y) & \text{if } x \neq y,
\end{cases}
\end{equation}
and most generalized contraction fixed-point theorems on $(X, \rho)$
follow from their classical analogues on $(X, \drho)$. The Matthews
theorem itself,
\begin{equation}\label{eq:matthews-banach}
\rho(f(x), f(y)) \le k\, \rho(x, y) \qquad (k \in [0, 1)),
\end{equation}
yields the Banach contraction $\drho(f(x), f(y)) \le k\, \drho(x, y)$,
from which existence and uniqueness of the fixed point follow. What the
partial setting preserves on its own is the supplementary property
$\rho(x^{\ast}, x^{\ast}) = 0$. The property is meaningful in
semantics; it does not constitute a new mathematical theorem.

The $\S14.6$ reduction is scoped to generalized contraction
hypotheses of the form $\rho(f(x), f(y)) \le \phi(\rho(x,y), \dots)$
on pairs. It is silent on theorems whose hypotheses or conclusions
mention the canonical order
\begin{equation}\label{eq:lerho}
x \lerho y \iff \rho(x, y) = \rho(x, x)
\end{equation}
attached to $\rho$. The relation $\lerho$ expresses that the
information carried by $y$ refines that carried by $x$ and has no
analogue on $(X, \drho)$: the equality $\drho(x,x) = 0$ for all $x$
collapses~\eqref{eq:lerho} to point equality. The fixed-point theorem
for operators monotone with respect to $\lerho$ stated
in~\cite[Thm.~14.5]{kirkShahzad2014}, which characterizes the fixed
point as $\sup_{\lerho}\{f^{n}(x_{0})\}$, therefore sits outside the
perimeter of the $\S14.6$ reduction by construction.

The present note formulates a path-averaged condition in the spirit of
Fabiano~\cite{fabiano2025paMetric, fabiano2027pa}, placed on the
\emph{self-distances} $\rho(f^{k} x_{0}, f^{k} x_{0})$ along a single
monotone orbit rather than on mutual distances between two orbits. The
condition is strictly weaker than the Matthews contraction. Its role in
the main theorem is narrow: existence of the fixed point, the supremum
characterization, and minimality above $x_{0}$ follow from Kleene's
theorem~\cite{kleene1938, davey2002} on a chain-complete poset under
Scott-continuity along the orbit, with no role for PA; PA is used only
to force the self-distance sequence to zero and so to secure the
strong-fixed-point property $\rho(x^{\ast}, x^{\ast}) = 0$.

The paper is organized as follows. Section~\ref{sec:prelims} collects
the definitions and auxiliary facts about partial metrics, the
canonical order, the Kleene fixed-point lemma, and 0-completeness that
enter the main theorem or locate its position in the Matthews literature.
Section~\ref{sec:main} states and proves the main result, isolating at
each step the use of structure that does not survive the $\drho$
reduction. Section~\ref{sec:examples} exhibits a worked example on the
lattice of bounded intervals and a counterexample separating our class
from the Matthews class. Section~\ref{sec:discussion} closes with a
split of the reduction-resistance claim into three sub-claims, a
comparison with the Matthews contraction, and four open questions.

\section{Preliminaries}
\label{sec:prelims}

We fix notation and recall the facts about partial metrics that enter
the proof of the main theorem. The content is standard but gathered
here with the specific ingredients the proof will use in
Section~\ref{sec:main}: the canonical order $\lerho$, the Kleene
fixed-point lemma on chain-complete posets, and the notion of
0-completeness that frames the perimeter of the $\drho$ reduction this
paper sidesteps.

A partial metric on a non-empty set $X$ is a map $\rho : X \times X \to
\RR^{+}$ satisfying, for all $x, y, z \in X$:
\begin{description}
\item[(PM1)] $x = y \iff \rho(x, x) = \rho(x, y) = \rho(y, y)$,
\item[(PM2)] $\rho(x, x) \le \rho(x, y)$,
\item[(PM3)] $\rho(x, y) = \rho(y, x)$,
\item[(PM4)] $\rho(x, z) \le \rho(x, y) + \rho(y, z) - \rho(y, y)$.
\end{description}
When $\rho(x, x) = 0$ for all $x$, inequality~(PM4) collapses to the
classical metric triangle inequality and $\rho$ is a metric. The
induced metric $\drho$ defined by~\eqref{eq:drho} makes $(X, \drho)$ a
standard metric space; it erases the self-distance information that
$\rho$ carries.

The canonical order $\lerho$ defined by~\eqref{eq:lerho} is a partial
order on $X$. Reflexivity reads $\rho(x, x) = \rho(x, x)$. For
antisymmetry, assume $x \lerho y$ and $y \lerho x$: then $\rho(x, y) =
\rho(x, x)$ and $\rho(y, x) = \rho(y, y)$, so by (PM3) $\rho(x, x) =
\rho(y, y)$, and (PM1) gives $x = y$. For transitivity, assume $x
\lerho y$ and $y \lerho z$: inequality (PM4) gives
\begin{align*}
\rho(x, z) &\le \rho(x, y) + \rho(y, z) - \rho(y, y) \\
           &= \rho(x, x) + \rho(y, y) - \rho(y, y) = \rho(x, x),
\end{align*}
and (PM2) gives the reverse inequality $\rho(x, x) \le \rho(x, z)$.

One further fact about $\lerho$ supplies the quantitative input the
proof in Section~\ref{sec:main} uses.

\begin{proposition}\label{prop:selfdistance-monotone}
For all $x, y \in X$ with $x \lerho y$, we have $\rho(y, y) \le \rho(x,
x)$.
\end{proposition}

\begin{proof}
From $x \lerho y$ we have $\rho(x, y) = \rho(x, x)$. Applying (PM2)
at $(y, x)$ and (PM3) gives $\rho(y, y) \le \rho(y, x) = \rho(x, y) =
\rho(x, x)$.
\end{proof}

So $\lerho$-increase and self-distance decrease travel together: the
more a point is refined in the $\lerho$ sense, the smaller its
self-distance. The proof of the main theorem will build on this fact
alone to force decay of the self-distance along a monotone orbit.

We now turn to the second ingredient: a fixed-point theorem for
monotone operators on a poset. Call a poset $(X, \lerho)$
\emph{chain-complete} if every $\lerho$-increasing chain $\{x_{n}\}_{n
\ge 0}$ admits a supremum in $X$. Call an operator $f : X \to X$
\emph{Scott-continuous} on such a chain if $f(\sup_{n} x_{n}) =
\sup_{n} f(x_{n})$. Our use of Scott-continuity will be restricted to
the orbit chains of $f$ itself.

\begin{lemma}[Kleene, form we use]\label{lem:kleene}
Let $(X, \lerho)$ be a chain-complete poset, $f : X \to X$ a
$\lerho$-monotone map Scott-continuous on its own orbits, and $x_{0}
\in X$ a point such that $x_{0} \lerho f(x_{0})$. Then $x^{\ast} :=
\sup_{\lerho}\{f^{n}(x_{0})\}_{n \ge 0}$ exists, is a fixed point of
$f$, and is the least fixed point of $f$ above $x_{0}$ in $\lerho$.
\end{lemma}

The proof follows~\cite{kleene1938, davey2002} and is not repeated
here. Three features of our Lemma~\ref{lem:kleene} are strictly weaker
than their counterparts in the PL-semantics tradition from which the
vocabulary derives. First, Scott-continuity is required only along the
orbit of $x_{0}$, not on all directed sets; the proof uses only one
chain, and the broader hypothesis would be unnecessary. Second,
chain-completeness is stated for $\omega$-chains, not for arbitrary
directed sets, so the lemma is an $\omega$-cpo statement, not a dcpo
one. Third, the iteration starts from an arbitrary post-fixed point
$x_{0} \lerho f(x_{0})$, not from a global bottom $\bot$; the ``least
fixed point'' concluded is least above $x_{0}$, not least globally.
These weakenings broaden the applicability of the theorem but remove
the Kleene-computability reading of the fixed point from the picture.

The third notion is relevant only to locate the perimeter of the
reduction that this paper sidesteps, and plays no role in the proof of
the main theorem. The partial metric topology admits a notion of
completeness tailored to sequences whose mutual distance decays to
zero, due to Romaguera~\cite{romaguera2010}.

\begin{definition}\label{def:0-cauchy}
A sequence $\{x_{n}\}_{n \ge 0}$ in a partial metric space $(X, \rho)$
is called \emph{0-Cauchy} if $\lim_{n, m \to \infty} \rho(x_{n}, x_{m})
= 0$.
\end{definition}

\begin{definition}\label{def:0-complete}
A partial metric space $(X, \rho)$ is called \emph{0-complete} if
every 0-Cauchy sequence converges to some $x \in X$ with $\rho(x, x) =
0$.
\end{definition}

Haghi, Rezapour and Shahzad~\cite{haghi2013} and the cited passage of
Kirk-Shahzad establish that $(X, \rho)$ is 0-complete if and only if
$(X, \drho)$ is complete. This equivalence is what makes transparent
the reduction of Matthews contractions to Banach contractions via
$\drho$. The result that concerns us (Theorem~\ref{thm:main}) does not
rely on 0-completeness; it deploys the partial order instead, which is
out of the reach of $\drho$.

We work throughout in the Matthews symmetric setting. K\"unzi, Pajoohesh
and Schellekens~\cite{kunziEtAl2006} drop the symmetry axiom (PM3) and
study partial quasi-metrics, in which the self-distance lower bound
splits into a left variant and a right variant. In that asymmetric
setting the canonical order itself bifurcates, which changes the shape
of a PA-Kleene theorem; the asymmetric case is not pursued here.

\section{The PA-Matthews-Kleene theorem}
\label{sec:main}

We now state the path-averaged condition and the theorem that follows
from it. In Fabiano's original setting on a complete metric or
$b$-metric space~\cite{fabiano2025paMetric, fabiano2027pa}, the PA
condition replaces the pointwise contraction $d(f(x), f(y)) \le k\,
d(x, y)$ by a bound averaged over $N$ consecutive iterates, permitting
individual contraction ratios above $1$ provided a block average
contracts by some $\alpha < 1$. Our variant places the inequality on
self-distances along a single monotone orbit. The condition and the
conclusion are both tied to a single $\lerho$-monotone orbit of the
operator; the proof reads self-distance decay off the PA inequality
and reads convergence off the chain supremum, in that order.

\begin{definition}\label{def:pa-matthews}
Let $(X, \rho)$ be a partial metric space and $\lerho$ its canonical
order. A map $f : X \to X$ is called \emph{PA-Matthews with parameters
$(\alpha, N)$} for $\alpha \in [0, 1)$ and $N \ge 1$ if there exists
$x_{0} \in X$ with $x_{0} \lerho f(x_{0})$ such that:
\begin{description}
\item[(M)] $f$ is $\lerho$-monotone,
\item[(PA)] the orbit $\{f^{n}(x_{0})\}_{n \ge 0}$ satisfies, for
every $j \ge 0$,
\begin{equation}\label{eq:pa-autodist}
\sum_{k = jN}^{(j+1)N - 1} \rho\bigl(f^{k + 1}(x_{0}), f^{k + 1}(x_{0})\bigr)
\le \alpha \sum_{k = jN}^{(j+1)N - 1} \rho\bigl(f^{k}(x_{0}), f^{k}(x_{0})\bigr).
\end{equation}
\end{description}
\end{definition}

Condition~\eqref{eq:pa-autodist} bears on \emph{self-distances}
$\rho(f^{k} x_{0}, f^{k} x_{0})$, not on mutual distances between two
orbits. This is the structural difference with the classical PA-Fabiano
condition on a metric space: without supposing two distinct orbits, the
inequality exerts pressure only on the self-distance of the chain of
interest, which is precisely the quantity measured by the partial
metric and erased by $\drho$.

\begin{theorem}[PA-Matthews-Kleene]\label{thm:main}
Let $(X, \rho)$ be a partial metric space, $\lerho$ its canonical
order, and assume that $(X, \lerho)$ is chain-complete. Let $f : X \to
X$ be a PA-Matthews map with parameters $(\alpha, N)$ starting at
$x_{0}$, Scott-continuous on the orbit $\{f^{n}(x_{0})\}_{n \ge 0}$.
Then:
\begin{enumerate}
\item the supremum $x^{\ast} := \sup_{\lerho}\{f^{n}(x_{0})\}_{n \ge
      0}$ exists in $X$;
\item $f(x^{\ast}) = x^{\ast}$;
\item $\rho(x^{\ast}, x^{\ast}) = 0$;
\item $x^{\ast}$ is the least fixed point of $f$ above $x_{0}$ in
      $\lerho$.
\end{enumerate}
\end{theorem}

\begin{proof}
Set $a_{k} := \rho(f^{k}(x_{0}), f^{k}(x_{0}))$. The assumption $x_{0}
\lerho f(x_{0})$ together with monotonicity (M) yields $f^{k}(x_{0})
\lerho f^{k + 1}(x_{0})$ for all $k$, so by
Proposition~\ref{prop:selfdistance-monotone} the sequence $\{a_{k}\}$
is decreasing in $\RR^{+}$:
\begin{equation}\label{eq:ak-decroissante}
a_{k + 1} \le a_{k} \qquad \text{for all } k \ge 0.
\end{equation}
Fix $j \ge 0$. Inequality (PA) at block $j$ reads, after reindexing,
\begin{equation}\label{eq:aN-bound}
\sum_{k = jN + 1}^{(j+1)N} a_{k} \le \alpha \sum_{k = jN}^{(j+1)N - 1} a_{k}.
\end{equation}
By the decrease~\eqref{eq:ak-decroissante}, $a_{(j+1)N} \le a_{k}$ for
all $k \in \{jN+1, \ldots, (j+1)N\}$ and $a_{k} \le a_{jN}$ for all $k
\in \{jN, \ldots, (j+1)N-1\}$, so
\begin{equation*}
N \, a_{(j+1)N} \le \sum_{k = jN + 1}^{(j+1)N} a_{k}
\le \alpha \sum_{k = jN}^{(j+1)N - 1} a_{k}
\le \alpha N \, a_{jN},
\end{equation*}
giving $a_{(j+1)N} \le \alpha \, a_{jN}$. By induction on $j$,
$a_{jN} \le \alpha^{j} a_{0}$. Monotone decrease within each block
then gives $a_{k} \le a_{\lfloor k/N \rfloor \cdot N}$ for all $k \ge
0$, hence
\begin{equation}\label{eq:ak-decay}
a_{k} \le \alpha^{\lfloor k/N \rfloor} a_{0}
\qquad \text{for all } k \ge 0,
\end{equation}
so $a_{k} \to 0$.

The chain $\{f^{n}(x_{0})\}_{n \ge 0}$ is $\lerho$-increasing by
monotonicity, and by the chain-completeness hypothesis, the supremum
$x^{\ast}$ exists. For every $k \ge 0$, $f^{k}(x_{0}) \lerho x^{\ast}$
gives $\rho(f^{k}(x_{0}), x^{\ast}) = \rho(f^{k}(x_{0}), f^{k}(x_{0}))
= a_{k}$, so~\eqref{eq:ak-decay} gives $\rho(f^{k}(x_{0}), x^{\ast})
\to 0$.

For property (3), (PM2) applied to $(x^{\ast}, f^{k}(x_{0}))$ gives
$\rho(x^{\ast}, x^{\ast}) \le \rho(x^{\ast}, f^{k}(x_{0})) = a_{k}$,
which tends to zero; hence $\rho(x^{\ast}, x^{\ast}) = 0$.

For property (2), Scott-continuity of $f$ on the orbit gives
\begin{equation*}
f(x^{\ast}) = f\bigl(\sup_{n} f^{n}(x_{0})\bigr) = \sup_{n}
f^{n + 1}(x_{0}) = \sup_{n \ge 1} f^{n}(x_{0}) = x^{\ast},
\end{equation*}
the penultimate equality holding because a chain and its shift by one
share the same supremum.

For property (4), let $y^{\ast}$ be a fixed point of $f$ with $x_{0}
\lerho y^{\ast}$. By monotonicity, $f^{n}(x_{0}) \lerho f^{n}(y^{\ast})
= y^{\ast}$ for all $n$, so $y^{\ast}$ is an upper bound of the chain
and $x^{\ast} \lerho y^{\ast}$. Hence $x^{\ast}$ is the least fixed
point above $x_{0}$.
\end{proof}

\begin{remark}[Scope of reduction-resistance]
\label{rem:nonreducible}
The four conclusions of Theorem~\ref{thm:main} have different
behavior under the $\drho$ reduction.

(A) \emph{Existence of a fixed point.} Under the theorem's hypotheses,
$f^{n}(x_{0}) \lerho f^{m}(x_{0})$ for $n \le m$ gives
$\rho(f^{n}(x_{0}), f^{m}(x_{0})) = a_{n}$, so~\eqref{eq:ak-decay}
makes the orbit 0-Cauchy in the sense of Def.~\ref{def:0-cauchy}. By
Haghi--Rezapour--Shahzad~\cite{haghi2013}, 0-completeness of
$(X, \rho)$ is equivalent to completeness of $(X, \drho)$; if
$(X, \drho)$ is complete, the orbit has a $\drho$-limit, which under
Scott-continuity along the orbit is a fixed point of $f$. So existence
reduces whenever $(X, \drho)$ is complete.

(B) \emph{Supremum characterization $x^{\ast} = \sup_{\lerho}\{f^{n}(x_{0})\}$,
and minimality above $x_{0}$ (conclusions (1) and (4)).} The order $\lerho$
is erased by $\drho$. Any statement of the fixed point in terms of
$\lerho$ or $\sup_{\lerho}$ has no $\drho$-analogue by construction;
this is a statement about which structure the theorem is written in,
not a theorem about non-reducibility.

(C) \emph{Strong-fixed-point property $\rho(x^{\ast}, x^{\ast}) = 0$
(conclusion (3)).} Under completeness of $(X, \drho)$, the 0-Cauchy
argument in (A) produces a limit satisfying $\rho = 0$ on the
diagonal by Haghi--Rezapour--Shahzad~\cite{haghi2013}, so (3) also
reduces when $(X, \drho)$ is complete.

Under the standing hypothesis of chain-completeness of $(X, \lerho)$
alone, conclusions (A) and (C) acquire a residual non-reductive
content: there is in principle a gap between chain-completeness of
$(X, \lerho)$ and 0-completeness of $(X, \rho)$ within which the
$\drho$-argument sketched in (A) is unavailable. We do not exhibit
a space in that gap; see the discussion in
Section~\ref{sec:discussion}. The mathematically unconditional
reduction-resistant conclusion of the theorem is (B).

The role of (PA) in all of this is confined to forcing $a_{k} \to 0$;
without (PA), Lemma~\ref{lem:kleene} alone gives the supremum, the
fixed point, and the minimality claim, but not the strong-fixed-point
property.
\end{remark}

\begin{remark}[Comparison with Theorem~14.5 of \cite{kirkShahzad2014}]
\label{rem:kirkshahzad145}
Theorem~14.5 of~\cite{kirkShahzad2014} treats the pointwise case $N =
1$ of Definition~\ref{def:pa-matthews}. Its hypothesis is the
pointwise self-distance contraction $\rho(f^{n}(x_{0}), f^{n}(x_{0}))
\le k \, \rho(f^{n-1}(x_{0}), f^{n-1}(x_{0}))$ for all $n$, together
with $\lerho$-monotonicity of $f$ and $x_{0} \lerho f(x_{0})$, under
0-completeness of $(X, \rho)$; its conclusion is that the orbit
supremum $x^{\ast}$ is the \emph{unique} fixed point of $f$ in the
up-set $\{z : x_{0} \lerho z\}$, is maximal in $(X, \lerho)$,
converges to the orbit in $\mathcal{T}(\rho^{s})$, and satisfies
$\rho(x^{\ast}, x^{\ast}) = 0$. Theorem~\ref{thm:main} extends this
result in two directions: (i) the pointwise contraction is replaced
by the block-averaged inequality~\eqref{eq:pa-autodist} with arbitrary
$N \ge 1$, in the spirit of Fabiano~\cite{fabiano2025paMetric,
fabiano2027pa}; (ii) 0-completeness of $(X, \rho)$ is replaced by
chain-completeness of $(X, \lerho)$, which is in principle a weaker
hypothesis (see question~(iv) of Section~\ref{sec:discussion}). The
cost is that our conclusion delivers only minimality in the up-set,
not uniqueness.

Neither extension is witnessed by the examples in this paper.
Examples~\ref{ex:intervals} and~\ref{ex:counter} are both 0-complete
and both satisfy the pointwise ($N = 1$) hypothesis, so both are
already covered by Theorem~14.5. What
Example~\ref{ex:counter} witnesses is the gap between the pointwise
self-distance contraction of Theorem~14.5 and the \emph{Matthews
contraction on all pairs}~\eqref{eq:matthews-banach} of~\cite[Thm.~14.1]{kirkShahzad2014}:
the latter fails on the pair $(a_{0}, b)$ while the former holds along
the orbit. Separating Theorem~\ref{thm:main} from Theorem~14.5 by a
concrete space --- one on which the block-averaged case with $N \ge
2$ holds but the pointwise case fails, or one in the
chain-complete-vs-0-complete gap of question~(iv) --- is left open.
\end{remark}

\section{Examples}
\label{sec:examples}

The next two examples show that the theorem is both applicable to
classical constructions and strictly more permissive than the Matthews
contraction. The first is a worked instance on the lattice of bounded
intervals that satisfies both our hypotheses and the Matthews
contraction, used to illustrate the supremum-based construction of the
fixed point. The second is a counterexample on a carefully constructed
partial metric space, satisfying our hypotheses but failing the
Matthews contraction on a specific pair of points, which establishes
that our class is strictly larger than the Matthews class on at least
one space.

We begin with the worked example. Let $\mathcal{I}$ be the family of
closed bounded intervals $[a, b] \subseteq \RR$ with $a \le b$,
equipped with the partial metric
\begin{equation}\label{eq:rho-intervals}
\rho([a, b], [c, d]) = \max(b, d) - \min(a, c),
\end{equation}
see Kirk-Shahzad~\cite[\S14.2.3]{kirkShahzad2014}. The self-distance
$\rho([a, b], [a, b]) = b - a$ is the length of the interval,
interpreted as incompleteness of the information: the singleton $[a,
a]$ carries zero self-distance. The canonical order $\lerho$ satisfies
$[a, b] \lerho [c, d]$ if and only if $[c, d] \subseteq [a, b]$, i.e.\
refinement by reverse inclusion.

\begin{example}\label{ex:intervals}
On the interval lattice $\mathcal{I}$ just introduced, define
$f : \mathcal{I} \to \mathcal{I}$ by $f([a, b]) = \left[a,
\tfrac{a + b}{2}\right]$. Starting at $x_{0} = [0, 1]$, induction gives
$f^{n}(x_{0}) = [0, 2^{-n}]$. Four quick verifications tie the example
to Theorem~\ref{thm:main}. First, $x_{0} = [0, 1] \lerho f(x_{0}) =
\left[0, \tfrac{1}{2}\right]$ since $\left[0, \tfrac{1}{2}\right]
\subseteq [0, 1]$. Second, $f$ is $\lerho$-monotone: if $[c, d]
\subseteq [a, b]$, then $c \ge a$ and $d \le b$, hence $f([c, d]) =
\left[c, \tfrac{c + d}{2}\right] \subseteq \left[a,
\tfrac{a + b}{2}\right] = f([a, b])$. Third, self-distances are $a_{n}
= 2^{-n}$ and satisfy the PA condition~\eqref{eq:pa-autodist} with $N
= 1$ and $\alpha = \tfrac{1}{2}$: $a_{1} = \tfrac{1}{2} = \tfrac{1}{2}
\cdot 1 = \tfrac{1}{2} a_{0}$. Fourth, the
chain $\{[0, 2^{-n}]\}_{n \ge 0}$ (decreasing for inclusion) has
intersection $\{0\}$ and $[0, 0] \subseteq [0, 2^{-n}]$ for all $n$, so
$x^{\ast} = \sup_{\lerho}\{f^{n}(x_{0})\} = [0, 0]$; this is a fixed
point with zero self-distance. For completeness, $(\mathcal{I},
\lerho)$ is chain-complete: a $\lerho$-increasing chain $\{[a_{n},
b_{n}]\}$ is a nested decreasing-for-inclusion sequence whose
intersection $[\sup_{n} a_{n}, \inf_{n} b_{n}]$ is a nonempty closed
bounded interval, which is the desired supremum.
\end{example}

\begin{figure}[ht]
\centering
\begin{tikzpicture}[scale=0.9]
\draw[->] (-0.3, 0) -- (5, 0) node[right] {$\RR$};
\foreach \n/\len/\y in {0/4.0/3.2, 1/2.0/2.4, 2/1.0/1.6, 3/0.5/0.8} {
  \draw[very thick, blue] (0, \y) -- (\len, \y);
  \draw[blue, fill=blue] (0, \y) circle (1.5pt);
  \draw[blue, fill=blue] (\len, \y) circle (1.5pt);
  \node[left, blue] at (-0.1, \y) {$[0, 2^{-\n}]$};
}
\draw[red, fill=red] (0, 0.15) circle (2pt);
\node[right, red] at (0.15, 0.15) {$x^{\ast} = [0, 0]$};
\draw[->, thick, gray!70!black, decorate, decoration={snake, amplitude=0.4mm}]
  (4.2, 3.2) -- (2.2, 2.4) -- (1.2, 1.6) -- (0.7, 0.8) -- (0.2, 0.18);
\node[right] at (4.2, 3.5) {\small chain $\{f^{n}(x_{0})\}$};
\node[right] at (4.2, 3.0) {\small refinement $\lerho$};
\end{tikzpicture}
\caption{The $\lerho$-increasing chain $\{f^{n}(x_{0}) = [0, 2^{-n}]\}$
of Example~\ref{ex:intervals} converges to its supremum $x^{\ast} = [0,
0]$. The self-distance $\rho(f^{n}(x_{0}), f^{n}(x_{0})) = 2^{-n}$
tends to zero, giving $\rho(x^{\ast}, x^{\ast}) = 0$.}
\label{fig:lattice}
\end{figure}
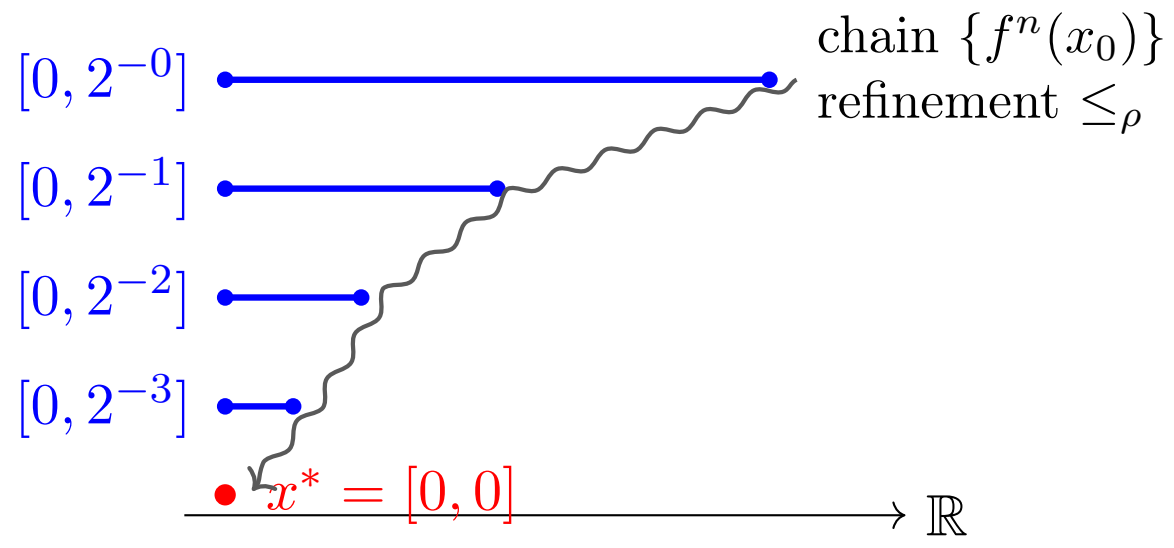

Along the chain $\{f^{n}(x_{0})\}$, Example~\ref{ex:intervals}
satisfies the Matthews contraction~\eqref{eq:matthews-banach} with $k
= \tfrac{1}{2}$; globally on $\mathcal{I}$ it does not, as witnessed
by $\rho(f([0,1]), f([1,2])) = \rho([0, \tfrac{1}{2}], [1, \tfrac{3}{2}])
= \tfrac{3}{2}$, which exceeds $\tfrac{1}{2} \, \rho([0,1], [1,2]) = 1$.
Either way the example does not separate the two classes. For separation
we need a space on which PA-Matthews succeeds while the Matthews
contraction fails on at least one pair of points. We construct one.

\begin{example}\label{ex:counter}
Let $X = \{a_{n} : n \ge 0\} \cup \{a_{\infty}\} \cup \{b\}$, where
the $a_{n}$, $a_{\infty}$, and $b$ are pairwise distinct formal
symbols, and define $\rho$ by
\begin{align*}
\rho(a_{n}, a_{n}) &= 2^{-n} &(n \ge 0), \\
\rho(a_{\infty}, a_{\infty}) &= 0, &\\
\rho(b, b) &= 10, &\\
\rho(a_{n}, a_{m}) &= 2^{-\min(n, m)} &(n \neq m\ \text{in } \NN), \\
\rho(a_{n}, a_{\infty}) &= \rho(a_{\infty}, a_{n}) = 2^{-n} &(n \ge 0), \\
\rho(a_{n}, b) &= \rho(b, a_{n}) = 10 &(n \ge 0), \\
\rho(a_{\infty}, b) &= \rho(b, a_{\infty}) = 10.
\end{align*}
Define $f : X \to X$ by $f(a_{n}) = a_{n + 1}$ for $n \in \NN$,
$f(a_{\infty}) = a_{\infty}$ and $f(b) = b$.
\end{example}

We verify six properties for this construction, all
straightforward but stated to make the counterexample self-contained.

First, $\rho$ is a partial metric. (PM1) and (PM3) are immediate.
(PM2) requires $\rho(x, x) \le \rho(x, y)$ for all $y$: for $x = a_{n},
y = a_{m}$, $\rho(a_{n}, a_{m}) = 2^{-\min(n, m)} \ge 2^{-n} =
\rho(a_{n}, a_{n})$; for $x = a_{n}, y = b$, $\rho(a_{n}, b) = 10 \ge
2^{-n} = \rho(a_{n}, a_{n})$; other cases are analogous. For (PM4) on
the triple $(a_{n}, a_{m}, a_{p})$ with $n \le m \le p$, LHS $= 2^{-n}$,
RHS $= 2^{-n} + 2^{-m} - 2^{-m} = 2^{-n}$. The other orderings and
cases involving $a_{\infty}$ or $b$ are checked similarly; for the
triple $(a_{n}, b, a_{m})$, LHS $= 2^{-\min(n, m)}$, RHS $= 10 + 10 -
10 = 10 \ge 2^{-\min(n, m)}$, which holds.

Second, the canonical order $\lerho$ has an explicit description. For
$n, m \in \NN$ with $n \le m$, $\rho(a_{n}, a_{m}) = 2^{-n} =
\rho(a_{n}, a_{n})$, hence $a_{n} \lerho a_{m}$. For $n \ge 0$,
$\rho(a_{n}, a_{\infty}) = 2^{-n} = \rho(a_{n}, a_{n})$, so $a_{n}
\lerho a_{\infty}$. For $(b, a_{n})$, $\rho(b, a_{n}) = 10 = \rho(b,
b)$, so $b \lerho a_{n}$; on the other hand $\rho(a_{n}, b) = 10 \neq
2^{-n} = \rho(a_{n}, a_{n})$, so $a_{n} \not\lerho b$. Thus $b$ sits
strictly below every $a_{n}$ and $a_{\infty}$. The poset has two
branches above $b$: an infinite chain $\{a_{n}\}$ rising to
$a_{\infty}$, and nothing else.

Third, $f$ is $\lerho$-monotone: $f(a_{n}) = a_{n+1} \lerho a_{m+1} =
f(a_{m})$ when $n \le m$; $f(b) = b \lerho a_{n+1} = f(a_{n})$ by the
above; $f(a_{\infty}) = a_{\infty}$ is maximal, automatically
compatible.

Fourth, the PA condition starting from $x_{0} = a_{0}$ holds with $N =
1$ and $\alpha = \tfrac{1}{2}$: self-distances along the orbit are
$a_{k} = 2^{-k}$, and $\rho(f(a_{0}), f(a_{0})) = 2^{-1} = \tfrac{1}{2}
\, \rho(a_{0}, a_{0})$. Also $a_{0} \lerho f(a_{0}) = a_{1}$ and $f$ is Scott-continuous
on $\{a_{n}\}_{n \ge 0}$: the chain reaches $a_{\infty}$ as its
supremum and $f(a_{\infty}) = a_{\infty}$.

Fifth, Theorem~\ref{thm:main} applies and gives $x^{\ast} = a_{\infty}$
as the fixed point with $\rho(a_{\infty}, a_{\infty}) = 0$.

Sixth, the Matthews contraction fails. For any $k \in [0, 1)$ take the
pair $(a_{0}, b)$:
\begin{equation*}
\rho(f(a_{0}), f(b)) = \rho(a_{1}, b) = 10 \not\le k \cdot 10 = k
\cdot \rho(a_{0}, b).
\end{equation*}
No $k < 1$ makes $f$ a Matthews contraction in the
sense~\eqref{eq:matthews-banach}, so Theorem~14.1 of
Kirk-Shahzad~\cite{kirkShahzad2014} does not apply to this $(X, \rho, f)$.

\begin{remark}[Scope of the separation]
The point $b$ plays the role of a ``parasite point'' sitting below the
PA chain in the $\lerho$ order but escaping any global Matthews
contraction. The PA-Matthews condition controls only the self-distance
along the chain and is insensitive to the presence of such parasite
points. The separation that Example~\ref{ex:counter} witnesses is
therefore \emph{off-chain}: PA succeeds on the orbit $\{a_{n}\}$ while
the Matthews contraction fails on the pair $(a_{0}, b)$, which lies
outside the orbit. The example does not exhibit a space on which the
Matthews contraction fails on the chain itself; whether such a space
can be constructed is left open.

Mild weakenings of the Matthews contraction in the literature reach
Example~\ref{ex:counter}. Berinde-type almost contractions on partial
metrics, $\rho(f(x), f(y)) \le k \, \rho(x, y) + L \, \rho(y, f(x))$,
hold on the pair $(a_{0}, b)$ with $k = 0, L = 1$ since $\rho(f(a_{0}),
f(b)) = \rho(a_{1}, b) = 10 = 0 \cdot \rho(a_{0}, b) + 1 \cdot \rho(b,
a_{1})$. What such weakenings do not deliver on their own is the
$\sup_{\lerho}$ characterization or the strong-fixed-point property,
which are the content specific to our hypothesis.
\end{remark}

\section{Discussion}
\label{sec:discussion}

We close with three observations: where the proof uses structure not
available on $(X, \drho)$ (reduction-resistance), how our class
compares with the Matthews class (strict containment at the price of
global uniqueness), and what remains open.

The first observation concerns the perimeter of the $\drho$ reduction
relative to Theorem~\ref{thm:main}. As
Remark~\ref{rem:nonreducible} articulates, the four conclusions of
Theorem~\ref{thm:main} split into three groups by reduction behavior:
(A) existence and (C) the strong-fixed-point property both reduce to
standard arguments on $(X, \drho)$ whenever $(X, \drho)$ is complete,
by Haghi--Rezapour--Shahzad~\cite{haghi2013} combined with the
0-Cauchy argument applied to the orbit; (B) the $\sup_{\lerho}$
characterization and the minimality claim above $x_{0}$ are
reduction-resistant by construction, since $\lerho$ has no
$\drho$-analogue.

(B) is the mathematically unconditional carrier of
reduction-resistance in this paper. The residual non-reductive content
of (A) and (C) lives in the gap between chain-completeness of
$(X, \lerho)$ and 0-completeness of $(X, \rho)$: a partial metric
space can be chain-complete without being 0-complete, in which case
the $\drho$-argument sketched above is unavailable and the theorem's
existence and strong-fixed-point conclusions become genuinely
non-reductive. Exhibiting such a space, and the behavior of
Theorem~\ref{thm:main} on it, is an open direction; we flag it here
and record it in question~(iv) below.

The second observation compares our class with the Matthews class.
Theorem~14.1 of Kirk-Shahzad~\cite{kirkShahzad2014} assumes
contraction~\eqref{eq:matthews-banach} on all pairs $(x, y)$ and
concludes the fixed point via the $\drho$ reduction. The present
Theorem~\ref{thm:main} assumes nothing on pairs outside the chain, only
an averaged decrease of self-distances along one monotone orbit.
Example~\ref{ex:counter} exhibits a space on which our class strictly
contains the Matthews class. The price paid is the absence of global
uniqueness of the fixed point (Example~\ref{ex:counter} has an
additional fixed point $b$ that our result does not reach): we
guarantee minimality above the starting point, not uniqueness on all of
$X$. In applications where only the smallest fixed point above a given
seed is needed, this cost is acceptable. In applications where global
uniqueness matters, the Matthews contraction remains the appropriate
hypothesis.

Four open questions close the paper. (i) Can the PA condition on
self-distances be strengthened to a PA condition on all pairs $(x, y)
\in X \times X$ that would recover global uniqueness without falling
back into the $\drho$ reduction? One avenue is to localize PA to
chain-neighbourhoods and use a gluing argument across them. (ii) If $f$
represents an approximate program transformation in the sense of
Geoffroy and Pistone~\cite{geoffroy2021}, does the PA-Matthews
condition admit an interpretation as a cumulative approximation cost
that decreases on average? (iii) Asymmetric extension in the partial
quasi-metric setting of K\"unzi, Pajoohesh and
Schellekens~\cite{kunziEtAl2006}. The two canonical orders induced by
axioms (1a) and (1b) give two parallel theorems whose structural
relationship (independence, equivalence, two-sided fusion) is open.
(iv) Does there exist a partial metric space $(X, \rho)$ that is
chain-complete in $(X, \lerho)$ but for which $(X, \drho)$ is not
complete, satisfying the hypotheses of Theorem~\ref{thm:main}? Such a
space would witness the non-reductive content of conclusions (A) and
(C) of Remark~\ref{rem:nonreducible}. A related question: is the
separation in Example~\ref{ex:counter} achievable on-chain, i.e., does
there exist a space on which PA-Matthews holds along
$\{f^{n}(x_{0})\}$ but the Matthews contraction fails on some pair
within the orbit itself?

\end{document}